\documentclass[12pt, reqno]{amsart}

\author[P.~Leonetti]{Paolo Leonetti}
\address{Department of Decision Sciences, Universit\`a ``Luigi Bocconi'', via Roentgen 1, 20136 Milan, Italy}
\email{leonetti.paolo@gmail.com}
\urladdr{https://sites.google.com/site/leonettipaolo/}
\keywords{Cesàro convergence; strong Cesàro convergence.}
\subjclass[2010]{Primary: 40D25. Secondary: 40A05, 40G15.}

\title{A Characterization of Cesàro Convergence}

\usepackage[T1]{fontenc}
\usepackage{amsmath}
\usepackage{amssymb}
\usepackage{amsthm}
\usepackage[left=3.5cm, right=3.5cm, paperheight=11.8in]{geometry}
\usepackage{hyperref}
\usepackage{fancyhdr}
\usepackage{enumitem}
\usepackage{comment}
\usepackage{nicefrac}
\usepackage{bm}
\usepackage{mathrsfs}
\usepackage{graphicx}
\usepackage[utf8]{inputenc}
\usepackage{cancel}
\usepackage{mathtools}
\usepackage{tikz}
\usetikzlibrary{cd,decorations.pathreplacing,matrix,arrows,positioning,automata,shapes,shadows,calc,fadings,decorations,snakes,through,intersections}
\usepackage{float}

\AtBeginDocument{%
   \def\MR#1{}
}

\newtheorem{thm}{Theorem}[section]
\theoremstyle{definition} 
\newtheorem{defi}[thm]{Definition}%[section]
\let\olddefi\defi
\renewcommand{\defi}{\olddefi\normalfont}
\newtheorem{question}{Question}%[section]
\let\oldquestion\question
\renewcommand{\question}{\oldquestion\normalfont}
\newtheorem{example}[thm]{Example}
\let\oldexample\example
\renewcommand{\example}{\oldexample\normalfont}
\newtheorem{rmk}[thm]{Remark}
\let\oldrmk\rmk
\renewcommand{\rmk}{\oldrmk\normalfont}

\hypersetup{
    pdftitle={A characterization of Cesàro convergence},
    pdfauthor={},
    pdfmenubar=false,
    pdffitwindow=true,
    pdfstartview=FitH,
    colorlinks=true,
    linkcolor=blue,
    citecolor=green,
    urlcolor=cyan
}

\providecommand{\MR}[1]{}

\providecommand{\MR}{\relax\ifhmode\unskip\space\fi MR }

\providecommand{\href}[2]{#2}

\begin{document}

\maketitle
\thispagestyle{empty}

\begin{abstract} 
We show that a real bounded sequence $(x_n)$ is Cesàro convergent to $\ell$ if and only if the sequence of averages with indices in $[\alpha^k,\alpha^{k+1})$ converges to $\ell$ for all $\alpha>1$. If, in addition, the sequence $(x_n)$ is nonnegative, then it is Cesàro convergent to $0$ if and only if the same condition holds for some $\alpha>1$.
\end{abstract}

\section{Introduction.}\label{sec:intro}

A real sequence $x=(x_n)_{n\ge 1}$ is said to be \emph{Cesàro convergent} to $\ell \in \mathbf{R}$ if 
$$%\begin{equation}\label{eq:cesarodef}
\lim_{n\to \infty}\frac{1}{n}\sum_{k=1}^nx_n=\ell.
$$%end{equation}
This is a weaker notion than ordinary convergence: indeed, it involves a kind of smoothing of the original sequence by computing its partial averages. 
The detailed theory of Cesàro convergence is discussed in Hardy’s classic textbook \cite{MR1188874}. 
The vector space $w_1$ of \emph{strongly Cesàro convergent} sequences $x$ (i.e., the sequences such that $\lim_{n\to \infty}\frac{1}{n}\sum_{i=1}^n|x_n-\ell|=0$ for some $\ell \in \mathbf{R}$), endowed with the norm
$$
%\textstyle 
\|x\|=\sup_{n\ge 1}\frac{1}{2^n}
\sum_{2^{n}\le k<2^{n+1}}|x_k|,
%\sum_{k=2^{n}}^{2^{n+1}-1}|x_k|,
$$
turns out to be a Banach space; see \cite{MR225044}. A characterization of strong Cesàro convergence for bounded sequences with statistical convergence can be found, for example, in \cite[Theorem 2.1]{MR954458}; see \cite{MR1006746} for extensions with summability methods.  

In a different context, it is known that a set $A$ of positive integers has \emph{asymptotic density} $0$, that is, $\lim_{n\to \infty}\frac{1}{n}\#(A\cap [1,n])=0$ (see, e.g., \cite{MR4054777}), if and only if 
$$
\lim_{n\to \infty}\frac{1}{2^n}\#(A\cap [2^n,2^{n+1}))=0;
$$
see \cite[Lemma 3.1]{MR3391516} and compare with the proof of \cite[Theorem 1.13.3(a)]{MR1711328}. 
Here, $\# S$ stands for the cardinality of a set $S$. 
Identifying $A$ with the sequence $(x_n)_{n\ge 1}$ such that $x_n=1$ if $n \in A$ and $x_n=0$ otherwise, 
as in \cite{MR759670}, 
the latter result can be rephrased as follows: a $\{0,1\}$-valued sequence $(x_n)_{n\ge 1}$ is Cesàro convergent to $0$ if and only if 
\begin{equation}\label{eq:condition01}
\lim_{n\to \infty}\frac{1}{2^n}
\sum_{2^{n}\le k<2^{n+1}}x_k=0.
%\sum_{k=2^{n}}^{2^{n+1}-1}x_k=0.
\end{equation}
This happens, for example, if $(x_n)_{n\ge 1}$ is the sequence associated with $A=\bigcup_{n\ge 1}\{2^n+1,\ldots,2^n+n\}$. 
First of all, we note that the condition that the sequence is $\{0,1\}$-valued cannot be omitted: indeed, the sequence $(x_n)_{n\ge 1}$ defined by 
\begin{equation}\label{eq:counterexmaple}
x_n=
\begin{cases}
\,\,\,\,\,1\,\,\,\,&\text{ if }2^k \le n < 3\cdot 2^{k-1} \text{ for some }k \ge 1\\
\,-1 & \text{ otherwise},
\end{cases}
\end{equation}
%$x_n=1$ for all $n\in \mathbf{N}$ such that there exists $k \in \mathbf{N}$ with $2^k \le n < 3\cdot 2^{k-1}$, and $x_n=-1$ otherwise, 
satisfies the limit \eqref{eq:condition01} and, on the other hand, $\limsup_{n\to \infty}\frac{1}{n}\sum_{k=1}^n x_k=\frac{1}{4}$ (computed along the subsequence $(3\cdot 2^k)_{k\ge 1}$). This example implies that, even for bounded sequences, \eqref{eq:condition01} need not be equivalent to 
Cesàro convergence. % \eqref{eq:cesarodef}.

However, replacing the base $2$ in condition \eqref{eq:condition01} with \emph{all} bases $\alpha>1$, gives a characterization of Cesàro convergence for bounded sequences, which is the content of our main result: 
\begin{thm}\label{prop:cesaro}
A bounded sequence $(x_n)$ 
is Cesàro convergent to $\ell$ 
if and only if
$$
%\lim_{n\to \infty}\frac{1}{n}\sum_{i=1}^n x_i=\ell \,\,\,\,\,\Longleftrightarrow\,\,\, \,\,\
%\forall \alpha>1, \quad 
\lim_{k\to \infty}\, \frac{1}{\alpha^{k+1}-\alpha^k}\sum_{\alpha^k\le i <\alpha^{k+1}}x_i=\ell
%\lim_{k\to \infty} \frac{1}{\#I_{\alpha, k}}\sum_{i \in I_{\alpha, k}}x_i=\ell
\quad 
\text{ for all }\alpha>1.
$$
\end{thm}

However, Theorem \ref{prop:cesaro} does not explain why only the base $\alpha=2$ appears in \eqref{eq:condition01}. In this regard, note that a \emph{nonnegative} sequence $(x_n)_{n\ge 1}$ is Cesàro convergent to $0$ if and only if it is strongly Cesàro convergent to $0$. Hence, condition \eqref{eq:condition01} is justified by the following:
\begin{thm}\label{prop:cesaro2strongly}
A bounded nonnegative sequence $(x_n)$ 
is Cesàro convergent to $0$ 
\textup{(}hence, strongly Cesàro convergent to $0$\textup{)} 
if and only if 
%there exists $\alpha>1$ such that 
$$
%\lim_{n\to \infty}\frac{1}{n}\sum_{i=1}^n x_i=\ell \,\,\,\,\,\Longleftrightarrow\,\,\, \,\,\
%\exists \alpha>1, \quad 
%\lim_{k\to \infty} \frac{1}{\#I_{\alpha, k}}\sum_{i \in I_{\alpha, k}}x_i=\ell.
\lim_{k\to \infty}\, \frac{1}{\alpha^{k+1}-\alpha^k}\sum_{\alpha^k\le i <\alpha^{k+1}}x_i=0 
\quad 
\text{ for some }\alpha>1.
$$
\end{thm}
Of course, Theorem \ref{prop:cesaro2strongly} holds by replacing $0$ with an arbitrary $\ell$ and using the hypothesis that $x_n\ge 0$ with $x_n\ge \ell$ for all $n$. We chose $\ell=0$ to ease the exposition. 

The proofs, which are completely elementary, follow in the next sections.

\section{Proof of Theorem \ref{prop:cesaro}.}

For each $\alpha>1$ and $j \in \mathbf{N}$, define 
$$
I_{\alpha, j}:=[\alpha^{j-1}, \alpha^j) \cap \mathbf{N}.
$$
Assume by convention that $\frac{1}{\# I_{\alpha,j}}:=0$ if $I_{\alpha, j}$ is empty. Thus, for each $\alpha>1$, we have $\lim_{j\to \infty}\frac{1}{\# I_{\alpha,j}}\cdot (\alpha^j-\alpha^{j-1})=1$. Moreover, for each $n,j \in \mathbf{N}^+$, set 
$$
a_n:=\frac{1}{n}\sum_{i=1}^n x_i\,\,\,\text{ and }\,\,\,b_j:=\frac{1}{\#I_{\alpha, j}}\sum_{i \in I_{\alpha, j}}x_i.
$$ 
Finally, note that we can suppose without loss of generality that $\ell=0$. 

\medskip

\textsc{If part.} Suppose that $\lim_{j\to \infty} b_j=0$. Fix $\varepsilon>0$, so that there exists $t_0 \in \mathbf{N}^+$ such that $|b_t|<\varepsilon$ for all $t\ge t_0$. 
Define $\theta:=\sup_n|x_n|$ and assume that $\theta>0$; 
otherwise the statement is trivially true. 
Also fix $\alpha>\max\{1, \frac{2\theta}{\varepsilon}\}$, and $n,k \in \mathbf{N}^+$ such that $n\ge \alpha^{t_0}\max\{1, \frac{\theta}{\varepsilon}, \frac{1}{\theta}\}$ and $n \in I_{\alpha,k}$. 
Then
\begin{equation}\label{eq:identitycesaro}
na_n=\sum_{i=1}^n x_i=\sum_{j=1}^{k-1}w_j b_j+\sum_{i \in I_{\alpha,k}\cap\, [1, n]}x_i, 
\end{equation}
where $w_j:=\#I_{\alpha,j}$ for each $j\in \mathbf{N}^+$. 
Note that $\sum_{j\le t-1}w_j \le \alpha^t$ for all $t \in \mathbf{N}^+$.  
Since $n\ge \alpha^{t_0}$, we have $k-1\ge t_0$, so that
$$
a_n=\frac{1}{n}\left(\sum_{j=1}^{t_0-1}w_j b_j+\sum_{j=t_0}^{k-1}w_j b_j+\sum_{i \in I_{\alpha,k}\cap\, [1, n]}x_i\right).
$$
Considering that $a_n$ is the average of all $b_j$ (with $j\le k-1$), each repeated $w_j$ times, and the remaining $x_i$ (with $i \in I_{\alpha,k}\cap\, [1, n]$), we obtain that
\begin{equation}\label{eq:mainchaininqueality}
\begin{split}
|a_n| &\le \frac{1}{n}\left(\sum_{j=1}^{t_0-1}w_j |b_j|+\sum_{j=t_0}^{k-1}w_j |b_j|+\sum_{i \in I_{\alpha,k}\cap\, [1, n]}|x_i|\right)\\
&\le \frac{1}{n}\left(\sum_{j=1}^{t_0-1}w_j \theta+\sum_{j=t_0}^{k-1}w_j \varepsilon +\sum_{i \in [n(1-1/\alpha), n]}\theta\right)\\
&\le \frac{1}{n}\left(\theta \alpha^{t_0}+n\varepsilon +\#[n(1-1/\alpha), n]\cdot \theta\right)
\\&
\le \frac{\theta\alpha^{t_0}}{n}+\varepsilon+\frac{\theta}{\alpha}+\frac{1}{n}
\le \frac{\theta\alpha^{t_0}}{n}+\varepsilon+\frac{2\theta}{\alpha}
% \\
%\le \frac{C\alpha^{t_0}}{n}+\varepsilon+\left(\frac{1}{\alpha}+\frac{\varepsilon}{2C}\right)\cdot C 
\le \varepsilon+\varepsilon+\varepsilon.
\end{split}
\end{equation}
%where $\|x\|=\sup_n |x_n|<\infty$ and $C:=\|x\|+|\ell|$. 
%Since $\alpha$ can be chosen such that $\varepsilon+\frac{\|x\|+|\ell|}{\alpha}\le \varepsilon$ and 
%Since there exists $n_0 \in \mathbf{N}^+$ such that $\frac{C\alpha^{t_0}}{n}\le \varepsilon$ for all $n\ge n_0$, it follows that $|a_n-\ell| \le 3\varepsilon$ for all $n\ge \max\{\alpha^{t_0}, n_0\}$. 
In the second line, we used that $I_{\alpha,k}\cap [1,n]$ is contained in $[\alpha^{k-1}, \alpha^k]\cap [1,n]$, which is, in turn, contained in $[n(1-\frac{1}{\alpha}),n]$; 
%since $\max I_{\alpha,k} \le \alpha\min I_{\alpha,k}$; 
also, in the last line, we used the inequalities $n\ge \alpha^{t_0}\frac{\theta}{\varepsilon}$, $\frac{1}{n} \le \frac{\theta}{\alpha^{t_0}} \le \frac{\theta}{\alpha}$, and $\alpha\ge \frac{2\theta}{\varepsilon}$. 

By the arbitrariness of $\varepsilon$, we conclude that $\lim_{n\to \infty} a_n=0$. 

\medskip

\textsc{Only If part.} Suppose that $\lim_{n\to \infty} a_n=0$ and fix $\alpha>1$. For each $j \in \mathbf{N}$, define $\iota_j:=\lceil \alpha^{j-1}\rceil$ and note that $\iota_j=\min I_{\alpha,j}$ if $j$ is sufficiently large.  Reasoning as in \eqref{eq:identitycesaro}, we obtain that
$$
(\iota_{k+1}-\iota_k)b_k=\sum_{i \in I_{\alpha,k}}x_i=\sum_{i < \iota_{k+1}}x_i-\sum_{i <\iota_k}x_i=a_{\iota_{k+1}-1}({\iota_{k+1}-1})-a_{\iota_{k}-1}({\iota_{k}-1})
$$
whenever $k$ is sufficiently large. This implies that
$$
\lim_{k\to \infty}b_k=\lim_{k\to \infty}a_{\iota_{k+1}-1}\cdot \lim_{k\to \infty}\frac{\iota_{k+1}-1}{\iota_{k+1}-\iota_k}-\lim_{k\to \infty}a_{\iota_{k}-1}\cdot \lim_{k\to \infty}\frac{\iota_{k}-1}{\iota_{k+1}-\iota_k}=0,
$$
completing the proof.
%\end{proof}

\section{Proof of Theorem \ref{prop:cesaro2strongly}.}

\textsc{If part.} The proof proceeds along the same lines as the previous result. 
Here, we let $\alpha>1$ be \emph{any} fixed number (in particular, not necessarily greater than $2\theta/\varepsilon$) and $n,k \in \mathbf{N}^+$ are still taken such that $n\ge \alpha^{t_0}\max\{1, \frac{\theta}{\varepsilon}, \frac{1}{\theta}\}$ and $n \in I_{\alpha,k}$. 
%\textcolor{red}{In 
%the only difference being in 
In addition, we have 
the following upper estimate: 
$$
\sum_{i \in I_{\alpha,k} \cap \,[1,n]}|x_i|
\le \sum_{i \in I_{\alpha,k}}x_i 
%=\sum_{\alpha^{k-1}\le i\le n} x_i
%= x_{\lceil \alpha^{k-1}\rceil}+\cdots+x_n
=b_k
%\le \varepsilon w_k 
\le \varepsilon (\alpha^k-\alpha^{k-1})
%\le \varepsilon \left(\alpha - \frac{1}{\alpha}\right) n
\le \varepsilon \alpha^k
\le \varepsilon \alpha n,
$$
where the first inequality depends on the fact the sequence has nonnegative terms. 
Therefore, 
%for any fixed $\alpha>1$, 
in place of the chain of inequalities \eqref{eq:mainchaininqueality},
we obtain that 
$$
|a_n| \le \frac{\theta \alpha^{t_0}}{n}+\varepsilon+\varepsilon
$$
which is smaller than $3\alpha\varepsilon$ if $n$ is sufficiently large. 

\medskip

\textsc{Only If part.} This follows by Theorem \ref{prop:cesaro}. % and the observation preceeding the statement of Theorem \ref{prop:cesaro2strongly}.

\section{Concluding Remarks.} 

In light of Theorem \ref{prop:cesaro2strongly}, one may hope for a characterization of the (upper or lower) asymptotic density of $A\subseteq \mathbf{N}$ in terms of (superior or inferior) limit of the block averages $\frac{1}{2^n}\#(A \cap [2^n,2^{n+1}))$, which was the original motivation for this work. However, this is not possible: the reason is along the same lines as the example given in \eqref{eq:counterexmaple}. 
Indeed, denote the upper and lower asymptotic density of a set $A\subseteq \mathbf{N}$ by
$$
\mathsf{d}^\star(A):=\limsup_{n\to \infty}\frac{1}{n}\#(A\cap [1,n])
\quad \text{ and }\quad 
\mathsf{d}_\star(A):=\liminf_{n\to \infty}\frac{1}{n}\#(A\cap [1,n]),
$$
respectively, see e.g. \cite{MR4054777}, and define 
$$
A_s:=\bigcup\nolimits_{n\ge 1}\{2^n+\lfloor s2^{n-1}\rfloor +t: t=1,\ldots,2^{n-1}\}
$$
for each $s \in [0,1]$. It follows that $\lim_n \frac{1}{2^n}\#(A_s \cap [2^n,2^{n+1}))=\frac{1}{2}$. On the other hand, none of the $A_s$ admits asymptotic density: in fact,
$$
\mathsf{d}_\star(A_s)=\frac{1}{2+s}<
\frac{2}{3+s}=\mathsf{d}^\star(A_s)
$$
% and $\mathsf{d}^\star(A_s)=\frac{2}{3+s}$ 
for each $s \in [0,1]$.  

\subsection{Acknowledgment.}
The author is grateful to the editor and the two anonymous referees for their remarks that allowed for a substantial improvement of the presentation.

%\nocite{*}
\bibliographystyle{amsplain}
%\bibliography{levy}

\end{document}